\documentclass[12pt,a4paper]{amsart}
\usepackage{amssymb,amsfonts,amsthm,amsmath}
\usepackage{graphicx} % for inclusion of image
\usepackage{color}

\theoremstyle{plain}

\numberwithin{equation}{section} \numberwithin{theorem}{section}
\numberwithin{lemma}{section} \numberwithin{definition}{section}
\numberwithin{corollary}{section} \textheight =24cm
\begin{document}
\title[Equations with digit reversal and powers]{Some equations with features of digit reversal and powers}
\author{Geoffrey B Campbell}
\address{Mathematical Sciences Institute \\
         The Australian National University \\
         ACT, 0200, Australia}
 \email{Geoffrey.Campbell@anu.edu.au}

\author{Aleksander Zujev}
\address{Physics Department \\
         University of California \\
         Davis, California 95616, USA}
\email{azujev@ucdavis.edu}
\keywords{Cubic and quartic equations, Counting solutions of Diophantine equations, Higher degree equations; Fermat's equation.}
\subjclass{Primary: 11D25; Secondary: 11D45, 11D41}

\begin{abstract}
 In this paper we consider integers in base 10 like $abc$, where $a$, $b$, $c$ are digits of the integer,
  such that $abc^2 - (abc \cdot cba) \; = \; \pm  n^2$, where $n$ is a positive integer,
  as well as equations $abc^2  - (abc \cdot cba) \; = \; \pm  n^3$,
  and $abc^3 - (abc \cdot cba) \; = \; \pm  n^2$
 We consider asymptotic density of solutions.
We also compare the results with ones with bases different from 10.
\end{abstract}
\maketitle

\section{Introduction} \label{S:intro}

The following are some examples of the Diophantine equation
\begin{equation} \label{E:1.1}
 abc^2 - (abc \cdot cba) = \pm  n^2
\end{equation}

%%\noindent
Examples of (\ref{E:1.1}):
\begin{eqnarray}
&& 528^2 - (528\cdot 825) = - 396^2 \nonumber \\
&& 539- (539\cdot 935) = - 462^2 \nonumber \\
&& 825- (825\cdot 528) = 495^2 \nonumber \\
&& 1296- (1296\cdot 6921) = - 2700^2 \nonumber
\end{eqnarray}
.
%%\noindent
There is a similar Diophantine equation to (\ref{E:1.1}), namely.
\begin{equation} \label{E:1.2}
 abc^2 - (abc\cdot cba) = \pm  n^3
\end{equation}
%%\noindent
Examples of (\ref{E:1.2}):
\begin{eqnarray}
&& 48^2 - (48\cdot 84) = - 12^3 \nonumber \\
&& 2744^2 - (2744\cdot 4472) = - 165^3 \nonumber \\
&& 5632^2 - (5632\cdot 2365) = 264^3 \nonumber \\
&& 7128^2 - (7128\cdot 8217) = - 198^3 \nonumber.
\end{eqnarray}

%%\noindent
There is a similar Diophantine equation to (\ref{E:1.1}) and (\ref{E:1.2}), namely
\begin{equation} \label{E:1.3}
 abc^3 - (abc\cdot cba) = \pm  n^2
\end{equation}
%%\noindent
Examples of (\ref{E:1.3}):
\begin{eqnarray}
&& 101^3 - (101\cdot 101) = 1010^2 \nonumber \\
&& 626^3 - (626\cdot 626) = 15650^2 \nonumber
\end{eqnarray}

\section{Main Results}

\subsection{Computational results}

\subsubsection{Equation (\ref{E:1.1})}

Solutions to (\ref{E:1.1}) under $530000$:
\begin{eqnarray*}
&& 528^2-528 \cdot 825=-396^2                \\
&& 539^2-539 \cdot 935=-462^2 \\
&& 825^2-825 \cdot 528=495^2 \nonumber \\
&& 1296^2-1296 \cdot 6921=-2700^2 \nonumber \\
&& 21296^2-21296 \cdot 69212=-31944^2 \nonumber \\
&& 35904^2-35904 \cdot 40953=-13464^2 \nonumber \\
&& 39204^2-39204 \cdot 40293=-6534^2 \nonumber \\
&& 51483^2-51483 \cdot 38415=25938^2 \nonumber \\
&& 83259^2-83259 \cdot 95238=-31581^2 \nonumber \\
&& 100793^2-100793 \cdot 397001=-172788^2 \nonumber \\
&& 120213^2-120213 \cdot 312021=-151848^2 \nonumber \\
&& 131043^2-131043 \cdot 340131=-165528^2 \nonumber \\
&& 184093^2-184093 \cdot 390481=-194922^2 \nonumber \\
&& 197516^2-197516 \cdot 615791=-287430^2 \nonumber \\
&& 214896^2-214896 \cdot 698412=-322344^2 \nonumber \\
&& 240426^2-240426 \cdot 624042=-303696^2 \nonumber \\
&& 243675^2-243675 \cdot 576342=-284715^2 \nonumber \\
&& 247192^2-247192 \cdot 291742=-104940^2 \nonumber \\
&& 251256^2-251256 \cdot 652152=-317376^2 \nonumber \\
&& 252486^2-252486 \cdot 684252=-330174^2 \nonumber \\
&& 262086^2-262086 \cdot 680262=-331056^2 \nonumber \\
&& 297992^2-297992 \cdot 299792=-23160^2 \nonumber \\
&& 324723^2-324723 \cdot 327423=-29610^2 \nonumber \\
&& 344619^2-344619 \cdot 916443=-443916^2 \nonumber \\
&& 348075^2-348075 \cdot 570843=-278460^2 \nonumber \\
&& 360639^2-360639 \cdot 936063=-455544^2 \nonumber \\
&& 371469^2-371469 \cdot 964173=-469224^2 \nonumber \\
&& 380208^2-380208 \cdot 802083=-400500^2 \nonumber \\
&& 382299^2-382299 \cdot 992283=-482904^2 \nonumber \\
&& 384659^2-384659 \cdot 956483=-468996^2 \nonumber \\
&& 395604^2-395604 \cdot 406593=-65934^2 \nonumber \\
&& 451737^2-451737 \cdot 737154=-359073^2 \nonumber \\
&& 456187^2-456187 \cdot 781654=-385323^2 \nonumber \\
&& 522729^2-522729 \cdot 927225=-459828^2 \nonumber \\
&& 523908^2-523908 \cdot 809325=-386694^2 \nonumber \\
&& 525625^2-525625 \cdot 526525=-21750^2 \nonumber \\
&& 528528^2-528528 \cdot 825825=-396396^2 \nonumber \\
\end{eqnarray*}

$\bullet$ Note that all three-digit numbers are multiples of 11:

$abc^2 - abc \cdot cba \; = \; 99 abc (a-c) \; = \; \pm  n^2$ must be multiple of $11^2$, so $abc$ must be multiple of $11$.

$\bullet$ Note the solutions 528 and 528528.
It is an instance of general rule:

If n-digit number $a$ is a solution of (\ref{E:1.1}), then so is $aa...a$, $k$ times concatenation of $a$.
It is simply proved, if consider that $aa...a$ = $((10^{nk}-1)/(10^n-1))$.
Considering that "n-digit number" may be interpreted widely - as have leading zeros,
then we have as solutions also
$a0a0...a$, $a00a00...a$, etc.

We therefore have an infinite number of solutions of (\ref{E:1.1}).

\subsubsection{Equation (\ref{E:1.2})}

Solutions to (\ref{E:1.2}) under $10^8$:
\begin{eqnarray}
&& 48^2-48 \cdot 84=-12^3 \nonumber \\
&& 2744^2-2744 \cdot 4472=-168^3 \nonumber \\
&& 4125^2-4125 \cdot 5214=-165^3 \nonumber \\
&& 5632^2-5632 \cdot 2365=264^3 \nonumber \\
&& 7128^2-7128 \cdot 8217=-198^3 \nonumber \\
&& 48000^2-48000 \cdot 84=1320^3 \nonumber \\
&& 49152^2-49152 \cdot 25194=1056^3 \nonumber \\
&& 148137^2-148137 \cdot 731841=-4422^3 \nonumber \\
&& 273273^2-273273 \cdot 372372=-3003^3 \nonumber \\
&& 321651^2-321651 \cdot 156123=3762^3 \nonumber \\
&& 456876^2-456876 \cdot 678654=-4662^3 \nonumber \\
&& 483153^2-483153 \cdot 351384=3993^3 \nonumber \\
&& 999000^2-999000 \cdot 999=9990^3 \nonumber \\
&& 3652264^2-3652264 \cdot 4622563=-15246^3 \nonumber \\
&& 5412825^2-5412825 \cdot 5282145=8910^3 \nonumber \\
&& 63936000^2-63936000 \cdot 63936=159840^3 \nonumber \\
\end{eqnarray}

$\bullet$ Note the solutions 999000.
It is part of a family of solutions 999000, 999999000000, ...,
or $10^{3k}(10^{3k}-1)$.
Therefore, there is an infinite number of solutions of (\ref{E:1.2}).
And this works for any base. For instance, in ternary system, equivalent solutions are
222000, 222222000000, ..., or, in base $b$, $b^{3k}(b^{3k}-1)$.

\subsubsection{Equation (\ref{E:1.3})}

Solutions to (\ref{E:1.3}) under $10^7$:
\begin{eqnarray}
&& 101^3-101 \cdot 101=1010^2 \nonumber \\
&& 626^3-626 \cdot 626=15650^2 \nonumber \\
&& 10001^3-10001 \cdot 10001=1000100^2 \nonumber \\
&& 1000001^3-1000001 \cdot 1000001=1000001000^2 \nonumber \\
&& 1040401^3-1040401 \cdot 1040401=1061209020^2 \nonumber \\
&& 2217122^3-2217122 \cdot 2217122=3301294658^2 \nonumber \\
&& 5053505^3-5053505 \cdot 5053505=11360279240^2 \nonumber \\
\end{eqnarray}
Notice that all solutions to (\ref{E:1.3}) are symmetric numbers.
It is easily explained. If number $m$ is symmetric, then

\noindent
$m^3 - m^2 = m^2(m-1) \; = \; n^2$, or $m-1 \; = \; k^2$.

\noindent
This has solutions more common, than when $m$ is not symmetric.

A family of symmetric solutions is $10^{2k} + 1$:
\begin{equation} \label{E:2.1}
(10^{2k} + 1)^3 - (10^{2k} + 1)^2 = ((10^{2k} + 1)10^k)^2
\end{equation}
Therefore, there is an infinite number of symmetric solutions of (\ref{E:1.3}).
%%The largest solution of (\ref{E:1.3}) that we found is 5053505, and then nothing at least below $10^8$.

\subsection{Results for bases different from 10}

We looked for solutions to (\ref{E:1.3}) in bases 3 - 9. We show results for bases 3 and 4.

Solutions in base 3 under $10^7$:
\begin{eqnarray*}
&& 101_{b3} \; 10^3-10 \cdot 10=30^2 \\
&& 222_{b3} \;  26^3-26 \cdot 26=130^2 \\
&& 10001_{b3} \;  82^3-82 \cdot 82=738^2 \\
&& 11202_{b3} \;  128^3-128 \cdot 184=1440^2 \\
&& 1000001_{b3} \;  730^3-730 \cdot 730=19710^2 \\
&& 2112112_{b3} \;  1850^3-1850 \cdot 1850=79550^2 \\
&& 100000001_{b3} \;  6562^3-6562 \cdot 6562=531522^2 \\
&& 101101101_{b3} \;  7570^3-7570 \cdot 7570=658590^2 \\
&& 222212222_{b3} \;  19601^3-19601 \cdot 19601=2744140^2 \\
&& 10000000001_{b3} \;  59050^3-59050 \cdot 59050=14349150^2 \\
&& 10112121101_{b3} \;  69697^3-69697 \cdot 69697=18400008^2 \\
&& 1000000000001_{b3} \;  531442^3-531442 \cdot 531442=387421218^2 \\
&& 2221000001222_{b3} \;  1555010^3-1555010 \cdot 1555010=1939097470^2 \\
&& 100000000000001_{b3} \;  4782970^3-4782970 \cdot 4782970=10460355390^2 \\
&& 100011000110001_{b3} \;  4862026^3-4862026 \cdot 4862026=10720767330^2 \\
&& 101102202201101_{b3} \;  5546026^3-5546026 \cdot 5546026=13060891230^2
\end{eqnarray*}

Solutions in base 4 under $10^7$:
\begin{eqnarray*}
&& 11_{b4} \;  5^3-5 \cdot 5=10^2 \\
&& 22_{b4} \;  10^3-10 \cdot 10=30^2 \\
&& 101_{b4} \;  17^3-17 \cdot 17=68^2 \\
&& 1001_{b4} \;  65^3-65 \cdot 65=520^2 \\
&& 2222_{b4} \;  170^3-170 \cdot 170=2210^2 \\
&& 10001_{b4} \;  257^3-257 \cdot 257=4112^2 \\
&& 11011_{b4} \;  325^3-325 \cdot 325=5850^2 \\
&& 20102_{b4} \;  530^3-530 \cdot 530=12190^2 \\
&& 100001_{b4} \;  1025^3-1025 \cdot 1025=32800^2 \\
&& 112211_{b4} \;  1445^3-1445 \cdot 1445=54910^2 \\
&& 202202_{b4} \;  2210^3-2210 \cdot 2210=103870^2 \\
&& 223322_{b4} \;  2810^3-2810 \cdot 2810=148930^2 \\
&& 1000001_{b4} \;  4097^3-4097 \cdot 4097=262208^2 \\
&& 10000001_{b4} \;  16385^3-16385 \cdot 16385=2097280^2 \\
&& 10100101_{b4} \;  17425^3-17425 \cdot 17425=2300100^2 \\
&& 100000001_{b4} \;  65537^3-65537 \cdot 65537=16777472^2 \\
&& 101202101_{b4} \;  71825^3-71825 \cdot 71825=19249100^2 \\
&& 110202011_{b4} \;  84101^3-84101 \cdot 84101=24389290^2 \\
&& 1000000001_{b4} \;  262145^3-262145 \cdot 262145=134218240^2 \\
&& 2212332122_{b4} \;  683930^3-683930 \cdot 683930=565610110^2 \\
&& 10000000001_{b4} \;  1048577^3-1048577 \cdot 1048577=1073742848^2 \\
&& 10010001001_{b4} \;  1065025^3-1065025 \cdot 1065025=1099105800^2 \\
&& 10122222101_{b4} \;  1157777^3-1157777 \cdot 1157777=1245768052^2 \\
&& 11002320011_{b4} \;  1322501^3-1322501 \cdot 1322501=1520876150^2 \\
&& 100000000001_{b4} \;  4194305^3-4194305 \cdot 4194305=8589936640^2 \\
&& 100120021001_{b4} \;  4293185^3-4293185 \cdot 4293185=8895479320^2 \\
&& 112120021211_{b4} \;  5866085^3-5866085 \cdot 5866085=14207657870^2 \\
&& 203122221302_{b4} \;  9284210^3-9284210 \cdot 9284210=28288987870^2
\end{eqnarray*}

We can see that almost all solutions are symmetric numbers. We discussed this for the base 10.
We here also see a family of symmetric solutions $10^{2k} + 1$, where $10$ is a base - 3 and 4 in decimal.
This solution works in any base.
Therefore, there is an infinite number of symmetric solutions of (\ref{E:1.3}) in any base.

\subsection{Asymptotic density of solutions}

We'll briefly discuss the density of the solutions for equations (\ref{E:1.1} - \ref{E:1.3}).

The first estimate of the density of the solutions for equation (\ref{E:1.1})
is to assume that $abc^2 - abc \cdot cba$ is a (quasi)random number.
The probability that the random number of the order $n^2$ is a perfect square is about $\frac{1}{2n}$.
Then the number of solutions under $n$ is approximately $\frac{1}{2} \log(n)$.
In fact, we have the number of solutions
under $10^5: \; 9$; under $10^6: \; 54$; under $10^7: \; 96$; under $10^8: \; 176$;
It may be approximately logarithmic dependence, but a few times more than predicted.

Similar estimate for the density of the solutions for equation (\ref{E:1.2})
gives for the density of solutions $\frac{1}{3n^{4/3}}$ and for
the number of solutions under $n$  $1 - \frac{1}{n^{1/3}}$.
Actually, there are 16 solutions under $10^8$.

And the estimate for the density of the solutions for equation (\ref{E:1.3})
gives for the density of solutions $\frac{1}{2n^{3/2}}$ and for
the number of solutions under $n$  $1 - \frac{1}{n^{1/2}}$.

These estimates are considerably lower than computational results.
This should mean that the numbers of type $abc^2 - abc \cdot cba$ and $abc^3 - abc \cdot cba$
are not entirely random.

As we discussed above, all three equations have special solutions, which makes the number of solutions infinite.
For proper estimate of the density of solutions, we should take into account these special solutions.
However, it is quite likely that we didn't find all special solutions.

As for the density of the solutions for equation (\ref{E:1.3}),
they are certainly not random - almost all solutions are symmetric numbers.
Among all $10^{2k}$ 2k-digit numbers, $10^k$ are symmetric numbers,
and among all $10^{2k+1}$ (2k+1)-digit numbers, $10^{k+1}$ are symmetric numbers.

%%A family of symmetric solutions is $10^{2k} + 1$:
%%\begin{equation} \label{E:2.1}
%%(10^{2k} + 1)^3 - (10^{2k} + 1)^2 = ((10^{2k} + 1)10^k)^2
%%\end{equation}
%%Therefore, there is an infinite number of symmetric solutions.
%%The largest solution of (\ref{E:1.3}) that we found is 5053505, and then nothing at least below $10^8$.
There are no non-symmetric solutions among those we found. Possibly our estimate of
the density of solutions, not counting symmetric ones, is correct in this case.

\end{document}